%% file: main.tex
\documentclass[12pt,a4paper]{article}
\usepackage[utf8]{inputenc}
\usepackage[english]{babel}
\usepackage[left=2.55cm, right=2.55cm, top=2.55cm, bottom=2.55cm]{geometry}
\usepackage{dsfont,amsmath,amssymb,stmaryrd}
\usepackage{delarray}
\usepackage{pgfplots}
\usepackage{graphicx}
\usepackage{csquotes}
\usepackage{multirow}
\usepackage{siunitx}
\usepackage[linesnumbered,ruled,vlined]{algorithm2e}
\usepackage{algpseudocode}
\usepackage[sorting=nyt, backend=biber, style=authoryear, natbib=true]{biblatex}

\let\cite\citet

\usepackage[hidelinks]{hyperref}

\DeclareFieldFormat{doi}{%
  \mkbibacro{DOI}\addcolon\space
  \ifhyperref
    {\href{https://doi.org/#1}{\nolinkurl{https://doi.org/#1}}}
    {\nolinkurl{#1}}}

\usepackage{cleveref}
\usepackage{todonotes}
\usepackage{authblk}
\crefname{inequality}{inequality}{inequalities}
\Crefname{inequality}{Inequality}{Inequalities}
\creflabelformat{inequality}{#2(#1)#3}
\hypersetup{
	pdftitle={A Decomposition Method for the Hybrid Quantum-Classical Solution of the Number Partitioning Problem},
	pdfauthor={Zongji Li, Tobias Seidel, Michael Bortz, Raoul Heese},
	pdfkeywords={Optimization, Quantum Annealing, number partitioning},
	pdfcreator={\LaTeX},
	pdfproducer={\LaTeX},
	pdfsubject={Optimization}
}
\bibliography{main}

\newcommand{\myCite}[2][1]{\cite{#2}}
\title{A Decomposition Method for the Hybrid Quantum-Classical Solution of the Number Partitioning Problem}
\author{Zongji Li , Tobias Seidel\thanks{Corresponding author. Email: tobias.seidel@itwm.fraunhofer.de}~, Michael Bortz , Raoul Heese}
\date{\today}
\affil{Fraunhofer Institute for Industrial Mathematics ITWM\\  67663 Kaiserslautern\\  Germany}

\begin{document}

\maketitle

\begin{abstract}
Current quantum computers can only solve optimization problems of a very limited size. For larger problems, decomposition methods are required in which the original problem is broken down into several smaller sub-problems. These are then solved on the quantum computer and their solutions are merged into a final solution for the original problem. Often, these decomposition methods do not take the specific problem structure into account. In this paper, we present a tailored method using a divide-and-conquer strategy to solve the number partitioning problem (NPP) with a large number of variables. The idea is to perform a specialized decomposition into smaller NPPs, which can be solved on a quantum computer, and then recombine the results into another small auxiliary NPP. Solving this auxiliary problem yields an approximate solution of the original larger problem. We experimentally verify that our method allows to solve NPPs with over a thousand variables using a D-Wave quantum annealer.
    
\end{abstract}

\input{body}

\printbibliography[heading=bibintoc]

\end{document}

%% file: body.tex
\section{Introduction} \label{sec:introduction}
We focus on the optimization version of the 2-way Number Partitioning Problem (NPP) (also known as the partition problem), one of the Karp's 21 NP-complete problems from \myCite{npcomplete}. Given a finite multiset 
\begin{equation} \label{eq:W}
W := \{w_1,\ldots, w_n\}
\end{equation}
of \(n \geq 1 \) elements \(w_i \in \mathbb{N}\) for all \(i \in \{1,\dots,n\}\), the 2-way NPP is defined as 
\begin{equation} \label{eq:NPP}
\min_{A \subseteq W} E(W, A)
\end{equation}
with the objective
\begin{equation} \label{eq:E}
E(W, A) := \left\vert \sum_{w \in A} w - \sum_{w' \in W \setminus A} w' \right\vert,
\end{equation}
which represents the error (or energy) of the solution \(A\) to the NPP generated by \(W\). Practical applications of NPP (or its extension to more than 2 sets known as multi-way NPP) include scheduling with multiple identical machines, as exemplified in \myCite{coffman} and \myCite{Pinedo.2012}.

As mentioned before, the NP-hardness of this problem has been established in \myCite{npcomplete}. Furthermore, it has been demonstrated that both exact cover and knapsack problems are linearly reducible (in terms of problem size) to NPPs without any constraints in \myCite[p. 100]{npcomplete}. This implies that solving the NPP also addresses these other NP-hard problems. Early exact algorithms with an exponential run-time are, for example, described in \myCite{horowitz} and \myCite{schroeppel}. 

Classical heuristics for NPP include methods such as the Largest Differencing Method (LDM) introduced by \myCite[p. 2]{kk}, and the Greedy algorithm (GR) as described in \myCite[p. 127]{mertens}. The LDM is considered to perform well on an average basis as a polynomial-time heuristic, as analyzed in \myCite{boettcher} and \myCite{michiels}. Based on LDM, \myCite{korf98} proposed an exact anytime algorithm for general \(l\)-way NPPs, which requires exponential time to obtain an optimal solution in the worst case. Alternatively, \myCite{krasecki2023} studies a hybrid approach using molecular computers.

In this study, we do not aim to solve the NPP on a classical device, but want to investigate to what extent quantum computers can find a solution. Solving optimization problems with the help of quantum computers is a highly relevant research topic with various opportunities and challenges. A more general discussion of this subject would go far beyond the scope of this paper. For a comprehensive review of this topic, we refer to \myCite{abbas2023quantum}.

To solve the NPP on a quantum computer, it has to be formulated in a suitable way. A popular approach is to define the problem as a Quadratic Unconstrained Binary Optimization (QUBO) model. A QUBO defines an optimization problem given by \begin{equation}\label{eq:qubo}
    \min_{x\in \{ 0, 1 \}^n} x^TQx,
\end{equation}
where \(Q \in \mathbb{R}^{n \times n}\) is a constant matrix. General properties about QUBOs and how to model problems as a QUBO can, for example, be found in \myCite{kochenberger2014} and \myCite{qubom}.

For the NPP, we can convert \Cref{eq:NPP} into a QUBO defined by \( Q \in \mathbb{Z}^{n\times n} \) with entries
\begin{equation}\label{eq:Q}
    Q_{ij}:=\begin{cases} 
    w_i w_j &\text{if }i\neq j \\ 
    w_i(w_i - c) &\text{if }i = j 
    \end{cases}
\end{equation}
for \(i,j \in \{1,\dots,n\}\) and \(c:=\sum_{w_i\in W}w_i\), as formulated, for example, in \myCite[p. 145-146]{qubom}. The \(n\)-dimensional vector \(x\) of binary optimization variables represents the solution \(A = \{w_i \in W \,\mid x_i = 1 \}\) of \Cref{eq:E}. The equivalence to \eqref{eq:NPP} can be seen by considering the squared error (Hamiltonian):
\begin{equation*}
    H:=\left(\sum_{i=1}^n w_i x_i - \sum_{i=1}^n w_i (1-x_i)\right)^2 = c^2 -4 c \sum_{i=1}^n w_i x_i +4 \left(\sum_{i=1}^n w_i x_i\right)^2 = c^2+x^\top Q x
\end{equation*}
and removing the constant part \(c^2\).

This QUBO representation enables us to solve an NPP on quantum devices or using Simulated Annealing on classical computers. For optimization on quantum devices, two major approaches exist. One using gate-based hardware and one using annealing-based hardware \myCite{alexeev2021}. In this study, we focus on the quantum annealer provided by \myCite{dwave2020}.

The current Advantage quantum annealer from D-Wave possesses over 5000 qubits, yet due to the limited connectivity not all of them can be utilized for solving a single problem instance. An additional step of finding a feasible minor embedding from the QUBO graph onto the hardware graph of qubits with limited connectivity is required, which is itself an NP-hard problem, as described in \myCite{zbinden}. Since the \(Q\) from \Cref{eq:Q} is fully connected (having no zero entries), finding a feasible embedding is in general hard for NPPs. Moreover, the resulting embedding is not guaranteed to be exact, potentially introducing errors into the solution even before the actual quantum annealing begins. Consequently, reducing problem size is essential in order to use these quantum machines for solving real-world problems, which can easily exceed this limit.

One hybrid decomposing solver \textit{qbsolv} was proposed by D-Wave itself, as detailed in \myCite{dwave1}. It is important to note that this general solver can be applied to all problems with an Ising formulation (and, equivalently, a QUBO model). Aimed at improving this method, the work in \myCite{okada} further exploits the connectivity within problems to accommodate larger sub-problems.

As also observed in \myCite{atobe}, \myCite{osaba}, \myCite{raymond23}, and \myCite{shay}, the current decomposition methods are designed for compatibility with as many distinct problem species as possible, although \myCite{bass} states that \textit{qbsolv} is in general still the best choice for large problems as of 2021.

In this paper, we propose an alternative to the current and general decomposition methods, a decomposition method that is specifically designed for the 2-way NPP. The main idea is to consider subsets of the multiset \(W\). These subsets form sub-NPPs. If they are chosen small enough, they can be solved on a current quantum device. Solving these sub-problems leads to partial solutions. In another step, these partial solutions are then used to formulate another auxiliary NPP. Subsequently, solving this auxiliary NPP allows us to merge the sub-solutions to an overall solution of the complete NPP. 

Numerical experiments show that the results are improved significantly compared to a general purpose decomposition method. This indicates that by focusing on specific problems and leveraging their unique properties, one can significantly enhance solution quality. 

The remaining paper is structured as follows. We introduce the proposed algorithm in \Cref{sec:method}. We then investigate the proposed method on a theoretical level in \Cref{sec:theory}, examining the influence of our decomposition method on so-called \emph{perfect solutions}. To further support our conclusions, we show the results of the numerical experiments in \Cref{sec:numexp}, comparing LDM, GR, Simulated Annealing, \textit{qbsolv}, and different variants of our algorithm. We conclude the paper in \Cref{sec:conclusions} with a summary of our results and ideas for future research directions.

\section{Method}\label{sec:method}
We consider the NPP given in \Cref{eq:NPP} with the error function \(E(W, A)\) and a given multiset \(W\) as introduced in \Cref{eq:E} and \Cref{eq:W}, respectively. The proposed method consists of three major steps: First: decomposition of the original NPP into multiple smaller sub-NPPs. Second: solving the sub-NPPs with Quantum Annealing (QA). Third and lastly: recombining the solutions of the sub-NPPs into a solution for the original NPP by solving an auxiliary NPP with Simulated Annealing (SA). This auxiliary NPP is in particular smaller than the original NPP. Both QA and SA utilize the formulation of a NPP as a QUBO specified in \Cref{eq:Q}. One of the core ideas is that solving the sub-NPPs in the second step and also the auxiliary NPP in the final step is easier than solving the original NPP. We describe the different steps in more detail in the following. To ease the notations, we define \([n] := \{1, 2, \ldots, n\}\) for a \(n\in \mathbb{N}\).

\textbf{First step (decomposition):}
For the first step, we consider for some \(m\in \mathbb{N}\) with \(m \leq n\) a decomposing vector \(j \in [m]^n\). This decomposing vector determines the sub-NPP as follows, we let:
\begin{equation*}
    W_k := \{w_i \in W \mid j_i = k\} \text{ for all } k \in [m]
\end{equation*}
This means that the  multiset \(W\) of the original NPP, \Cref{eq:W}, is decomposed into \(m \) smaller multisets \(\mathcal{W} := \{ W_1, \ldots , W_m \}\) with \(\bigcup_{k=1}^m W_k = W\) and \(\sum_{k=1}^m |W_k| = |W|\). Without loss of generality, we presume that \(\vert W_k \vert \geq 1\) for all \(k \in [m]\). Each of the multisets in \(\mathcal{W}\) can be used to construct a separate NPP that is smaller than the original NPP.

There are different possible choices on how to construct a decomposing vector \(j \in [m]^n\). In general, we would like to keep the size of each sub-problem equally large. This way, all the resulting sub-NPPs will be, for a suitably chosen \(m\), small enough to be solved on a quantum device in the following step.

\textbf{Second step (solving sub-problems):}
For each of the multisets \(W_k, k \in [m]\) we can now consider a separate NPP:
\begin{equation*}
    \min_{A\subseteq W_k} E(W_k,A).
\end{equation*}
We can now independently solve each NPP and obtain a solution \(A_k\). This NPP contains less elements than the original NPP. For an appropriately chosen decomposing vector \(j\), we use QA to obtain the sub-solutions. The resulting partitions---which are not necessarily the optimal partitions---are denoted by \(A_k^1 := A_k\) and \(A_k^2 := W_k \setminus A_k\), respectively, for all \(k \in [m]\). Without loss of generality, we presume \(\sum_{w\in A_k^1} w \geq \sum_{w'\in A_k^2} w'\).
We write \(\mathcal{A} := \{A_1^1,\dots,A_m^1\}\) to denote the full set of partitions.

\textbf{Third step (merging):}
Finally, we merge these partial sub-solutions to retrieve a solution for the original NPP. After the second step, we have a partitioning for each of the sub-NPPs. A naive way to combine them would be to simply take the unions \(\bigcup_{k=1}^m A_k^1\) and  \(\bigcup_{k=1}^m A_k^2\). However, this way, the errors accumulate. To combine the sets in such a way that errors cancel out, we are led to another number partitioning problem using the errors of the solutions of the sub-NPPs.

Let \(E_k:= |\sum_{w \in A_k^1}w - \sum_{w \in A_k^2}w| =\sum_{w \in A_k^1}w - \sum_{w \in A_k^2}w\), where we may omit the absolute value as we assumed the first sum to be larger. We construct another multiset \(W' := \{ E_k \,\mid k\in [m] \}\) and obtain the corresponding auxiliary NPP:
\begin{equation*}
    \min_{A\subseteq W'} E(W',A)
\end{equation*}
We utilize SA for this computation. While the resulting partition \(A'\) is again not necessarily the optimal solution, we typically expect that this NPP is easier to solve than the original NPP. First, it contains fewer elements. Second, the values of the elements are typically smaller, as they are the error values of another NPP.

For the actual solution $\hat{A}$ of the original NPP given in \Cref{eq:NPP}, we use the following reconstruction:
\begin{equation}\label{eq:Aest}
    \hat{A} := \hat{A}(\mathcal{A}, W', A') := \bigcup_{k=1}^{m} A_k^{u(E_k,A')}
\end{equation}
with the abbreviation
\begin{equation}
    u(E_k,A') := \begin{cases} 1 & \text{if } E_k \in A' \\ 2 & \text{if } E_k \notin A'\end{cases}
\end{equation}
as an estimate for the optimal partition \(A\) of the original NPP based on \(\mathcal{A}\) and \(A'\). The corresponding objective reads \(E' = \vert \sum_{k=1}^{m} E_k (-1)^{u(E_k,A')} \vert = E(W, \hat{A})\) by construction. However, it is again not guaranteed that \(\hat{A}\) is an optimal solution. The proposed method is summarized in \Cref{alg:partition}.

\begin{algorithm}[h!]
    \caption{Proposed algorithm to solve the NPP}\label{alg:partition}
    \DontPrintSemicolon
    \KwData{\(W, j\)}
    \KwResult{partition \(A\)}
    Initialize an empty list \(\mathcal{A}\)\;
    Initialize an empty list \(W'\)\;
    \(\mathcal{W} \leftarrow\) Decompose \(W\) according to \(j\)\;
    \For{each \(W_k\) in \(\mathcal{W}\)}{
        \(A_k^1 \leftarrow\) Solve the NPP for \(W_k\) using QA\;
        Append \(A_k^1\) to \(\mathcal{A}\)\;
        \(E_k \leftarrow\) Compute the energy of the solution using \Cref{eq:E}\;
        Append \(E_k\) to \(W'\)\;
    }
    \(A' \leftarrow\) Solve the NPP for \(W'\) using SA\;
    \(\hat{A} \leftarrow\) Evaluate \Cref{eq:Aest} using \(\mathcal{A}\), \(W'\) and \(A'\)\;
    \Return \(\hat{A}\)\;
\end{algorithm} 

We comment briefly on two aspects of the algorithm. First, in line 5  we use QA, while in line 9 we use SA, although other choices are possible. For example, we could replace QA in line 5 by SA, or SA in line 9 by QA. We will investigate and compare the different variants using numerical experiments in \Cref{sec:numexp}.

Second, the proposed method differs in a key aspect from many other hybrid algorithms, such as \textit{qbsolv}. Other hybrid algorithms often run a loop in which one tries to improve the solution step by step. This process terminates when no further progress is obtained or when the maximum number of iterations is reached. In \Cref{alg:partition}, we partition the problem into sub-problems only once, solve them all, and then combine the solutions in a single step.

Before we turn to the theoretical and numerical analysis of our algorithm, we provide a short example illustrating our core idea. Consider the NPP given by
\[W := \{1, 1, 3, 4, 5, 6\} \text{ and a decomposing vector } j = (1,1,1,2,2,2)^{\intercal}\] with \(m=2\). In the first step, we perform a decomposition into the multisets \(W_1 := \{1, 1, 3\}\) and \(W_2 := \{4, 5, 6\}\) according to \(j\) and solve an NPP for each of them using QA in the second step. Let us assume this approach yields the partitions
\[A_1^1 := \{3\}, A_1^2 := \{1,1\}\text{ for }W_1\text{ as well as } A_2^1 := \{4, 5\}, A_2^2 := \{6\}\text{ for }W_2\]
with the objectives \(E_1 = 3 - 1 - 1 = 1\) and \(E_2 = 4 + 5 - 6 = 3\), respectively. That is, we arrive at the set of sub-partitions \(\mathcal{A}=\{A_1^1, A_2^1\}\). There are two distinct ways to merge these sub-partitions to obtain \(A\) for the original NPP, either 
\[A_1 := A_1^1 \cup A_2^1 = \{3, 4, 5\} \text{ or } A_2 := A_1^1 \cup A_2^2 = \{3, 6\}\]
For the first choice, we find \[E(W, A_1) = \vert 3 + 4 + 5 - 1 - 1 - 6 \vert = 4 = E_2 + E_1,\] and for the second, \[E(W, A_2) = \vert 3 + 6 - 1 - 1 - 4 - 5 \vert = 2 = E_2 - E_1.\] 
While for $A_1$ the errors simply add up, we can use $E_2$ in the second possibility to make the overall error smaller.

These outcomes align with the potential results of \(W' = \{E_1, E_2\}\) when treated as a NPP. We consequently solve this auxiliary NPP in the third step using SA. Let us assume that we obtain the solution \(A' = \{E_2\}\) with the objective \(E' = E(W, A_2)\). Hence, we find
\(\hat{A} = A_1^2 \cup A_2^1 = \{1, 1, 4, 5\}\) according to \Cref{eq:Aest} as an estimate for the partition \(A\) of the original NPP with the objective \(E(W, \hat{A}) = E' = E_2 - E_1 = 2\). This is in fact an optimal solution.
    
\section{Conservation of Perfect Solutions} \label{sec:theory}
In line 3 of \Cref{alg:partition}, we decompose the NPP into smaller sub-NPPs. An important question is how many sub-problems should be chosen. Before studying this question from a numerical perspective in \Cref{sec:numexp}, we want to give some insights from a theoretical point of view. To do so, we study the \emph{perfect solutions}. 

We differentiate between a \textquotedblleft perfect\textquotedblright{} solution and an \textquotedblleft optimal\textquotedblright{} solution. A solution \(A \subseteq W\) is called perfect, if and only if \(E(W, A) \in \{0, 1\}\). Whereas \(A\) is optimal, if and only if no \(A'\) exists such that \(E(W,A') < E(W, A)\). A perfect solution is always optimal, but not necessarily vice versa. While a perfect solution does not exist for all instances of NPP, an optimal one always does.

It is important to note that the proposed algorithm does not necessarily require perfect or optimal solutions of the sub-problems to perform effectively, as we see later in \Cref{sec:numexp}. The elements in \(W'\) merely need to be close to each other. However, analyzing the behavior of QA and SA is challenging since they are both probabilistic heuristics. Results from QA can additionally be influenced by errors resulting from inadequate quantum hardware. To simplify the theoretical analysis of our decomposition method, we assume that both QA and SA consistently yield optimal solutions.

Even under the assumption of the optimal solution, it is possible that we lose optimal or perfect solutions due to the use of decomposing vector \(j\). For example, let \(W(n) := \{2^t\mid t\in [n]\} \). We observe that the optimal solution is always \(E(W(n), \{2^n\}) = 1\), since \(2^n - 1 = \frac{2^n - 1}{2 - 1} = \sum_{t \in [n - 1]}2^t\). However, if we consider any decomposing vector \(j\) with \(m\geq 2\) and at least two elements in each subset \(W_k\), we always obtain a sub-optimal solution. 

This means that it is easily possible to construct problems for which \Cref{alg:partition} results in a sub-optimal solution even under perfect conditions. To avoid this worst case scenario, we turn to a probabilistic analysis. We make the following simplifying assumptions:
\begin{itemize}
    \item The elements of \(W\) are drawn from a uniform distribution \(U[\lambda]\) in \([\lambda]\) for a \(\lambda \in \mathbb{N}\)
    \item The sum of elements is even (thus, a perfect solution will have \(E=0\))
\end{itemize}
The expected number of perfect solutions under these assumptions, as in \myCite[p. 436-437]{analysis}, is given by 
\begin{equation}\label{eq:sol}
    |Sol(W)| = \frac{2^n}{\lambda}
\end{equation} for a sufficiently large \(n = |W|\) and \(\lambda\) the maximum of the uniform distribution.

For simplicity in the notations, let us further assume that \(m \mid n\), such that the size of each sub-problem will be \(|W_k| = \frac{n}{m}\). As each element of \(W\) is assumed to be drawn from a uniform distribution, the same holds true for the subset $W_k$ and we can apply \Cref{eq:sol} on \(W_k\) (if we assume \(|W_k|\) is even). Consequently, for each \(|W_k|\) the expected number of perfect solutions is \(|Sol(W_k)| = \frac{2^{n/m}}{\lambda}\), which we transform into
\begin{equation}\label{eq:subsol}
    m = \frac{n}{\log_2(\lambda \times |Sol(W_k)|)}
\end{equation}

The expected number of solutions only appears in the logarithm. This means that, provided \(\lambda\) does not grow exponentially with \(n\), we can decompose the problem into many sub-problems and still expect many perfect solutions.

For example, if we set \(|Sol(W_k)| \geq 1\) as a condition, then the equation turns into
\begin{equation}\label[inequality]{ineq:finsol}
    m \leq \frac{n}{\log_2(\lambda)}
\end{equation}
Then, for \(n = 500\) and \(\lambda = 5000\), we compute \(m \leq 40\). 

Similarly, we can also evaluate the auxiliary NPP considered in the merging step. There we have two effects. First, the solution of the sub-NPPs makes the values for the merging smaller. This effect is more influential, according to \Cref{eq:sol}, if we choose \(m\) smaller and obtain better solutions for the sub-NPPs. On the other hand, by choosing \(m\) larger and thus having more elements in the NPP considered in the merging step, we also increase the expected number of perfect solutions in this step. This means that we have two competing aspects: we should choose a small \(m\) for good solutions of the sub-NPPs, and a large \(m\) for good solutions of the auxiliary NPP in the merging step.  However, as the expected number of perfect solutions grows exponentially in the problem size, we have, for large instances, a lot of flexibility in choosing \(m\), naturally under the assumption that \(\lambda\) does not grow exponentially relative to \(n\).

We observe that, according to \Cref{ineq:finsol}, our decomposition method begins to fail as \(\lambda \rightarrow 2^n\), but is expected to perform well with \(\lambda \in O(n^d)\) for a constant \(d \in \mathbb{N}\) and sufficiently large \(n\). 

However, \Cref{ineq:finsol} describes a worst-case scenario. We expect a perfect solution for the random instances with the \(m\) specified in \Cref{ineq:finsol}. It is indeed not a necessary condition for the result of \Cref{alg:partition} to be optimal or even perfect, since some of the assumptions made to compute this \(m\) are stringent. For example, \(W' = \{5, 6\}\) does not fulfill the assumption of \(w' \in W' \Rightarrow w' \in \{0, 1\}\), nonetheless, the resulting \(E\) corresponds to a perfect solution. We note that these cases are more difficult to analyze, due to the lack of a reasonable estimation of \(E\) of imperfect yet optimal solutions. Therefore, this study does not theoretically analyze the performance of \Cref{alg:partition} in complex NPP instances where \(\lambda \in o(2^n)\) and thus \(|Sol(W_k)|=0\).

\section{Numerical Experiments}\label{sec:numexp}
After these theoretical insights into our decomposition, we proceed with a numerical study. First, we describe how we generated the test instances and specify the implementation details. We then proceed with three different questions. In the first two experiments, we want to study two aspects of the decomposition: firstly, we study the influence of the sub-problem size on the solution quality; secondly, we examine the influence of SA and QA on the solutions of the different NPPs. In a final experiment, we study the main question of this work. We compare the solution quality of our problem-oriented decomposition method with the general-purpose decomposition implemented in \textit{qbsolv}. Along with this experiment, we also compare the results to the classical heuristics (LDM, GR, and SA) for solving the NPP. 

The considered instances are constructed as follows. We fix a number of elements \(n\in \mathbb{N}\). We then draw \(n\) elements at uniformly random from the interval \([5n, 10n]\). While the previous analysis suggests that an exponential growth might lead to harder problems, we note that the values for the matrix \(Q\) as given in \Cref{eq:Q} are the square of the values of the NPP. As these values increase, the precision of their mapping onto the quantum device declines.

We note that OUR, QB, and SA are probabilistic. Therefore, we run each of them 5 times to solve a single problem and include all 5 results in our evaluations (i.e., plots and tables). For each \(n\), we generate 10 problems. Thus, for probabilistic heuristics, we obtain a total of 50 data points for each \(n\) for the box plots. For deterministic LDM and GR, there are 10 data points for each \(n\).

In our experiments, we decompose the original NPP at random. To ensure a balanced size among the sub-NPPs, we randomly and uniformly shuffle the list derived from \(W\) and sequentially distribute its elements into \(m\) sub-lists, striving to keep the sizes of these sub-lists as similar as possible. Consequently, the sub-problems are derived from these sub-lists,

In the following, SA denotes Simulated Annealing as implemented in \myCite{sa}. QB refers to \textit{qbsolv}, as implemented in the dwave-hybrid python package, specifically under the name \textit{SimplifiedQbsolv}, using non-default parameters \(max\_iter=1\) and \(max\_subproblem\_size=40\). LDM and GR are the Largest Differencing Method and Greedy heuristics as described in \myCite[p. 2]{kk} and \myCite[p. 127]{mertens} (We implemented these algorithms using Python). And for our proposed \Cref{alg:partition}, we write OUR. No additional fine-tuning of parameters (including, but not limited to, chain strength and manual embedding) nor any pre- or post-processing techniques (such as local search) are applied in our experiments. 

For the first experiment, we investigate the influence of sub-NPP size on the quality of solutions represented by the error \(E\). We test OUR (\Cref{alg:partition}) with \(m \in \{\lfloor n/20 \rfloor, \lfloor n/40 \rfloor, \lfloor n/60 \rfloor, \lfloor n/80 \rfloor\}\), with OUR-N representing \(m = \lfloor n/N \rfloor\) respectively.

\begin{figure}[hbt!]
    \centering
    \includegraphics[width=\textwidth]{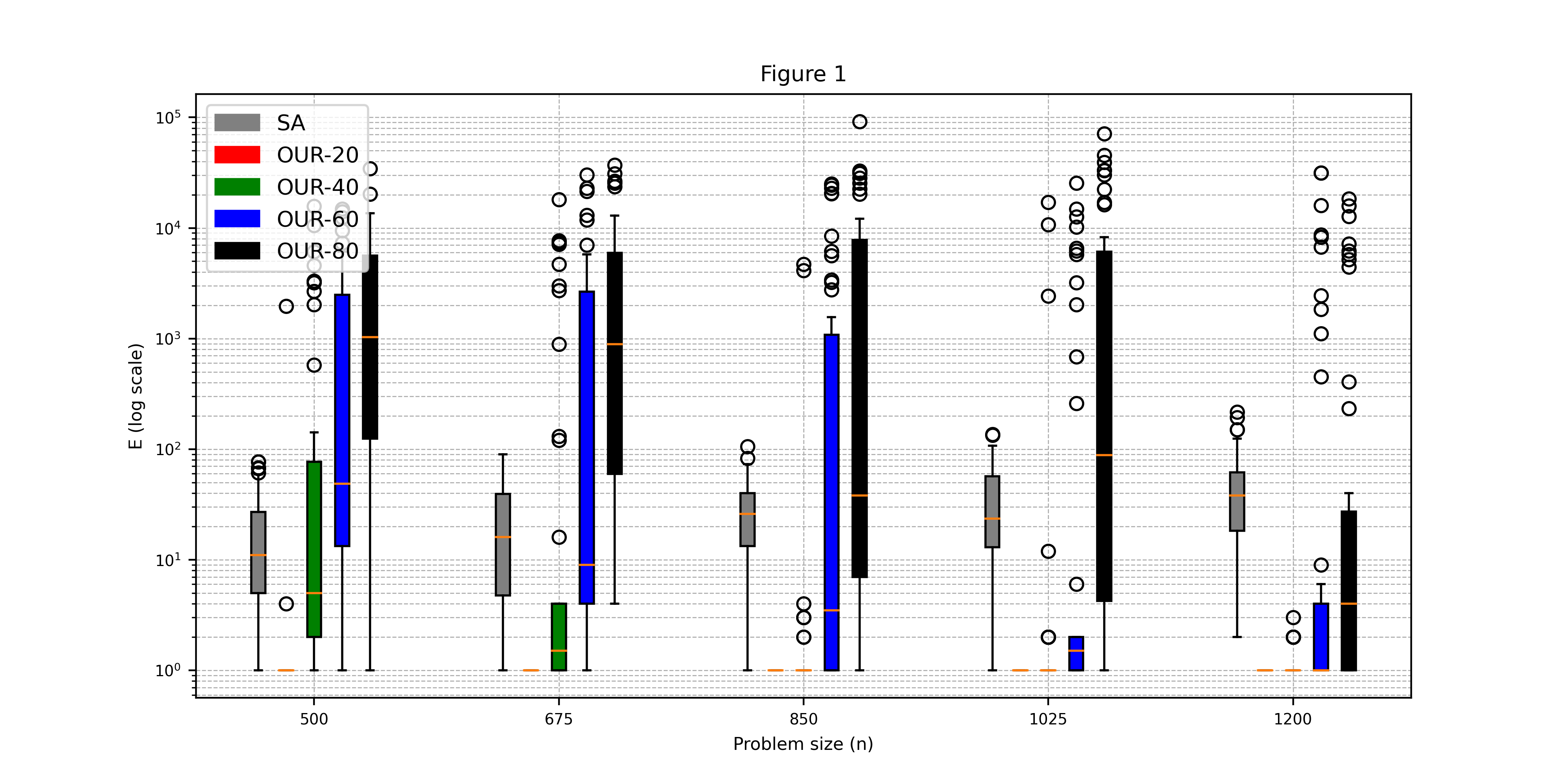}
    \caption{Experiment 1: comparing solution error \(E\) with different \(m\) for 5 different problem sizes}
    \label{fig:1}
\end{figure}

The experiments show that with OUR-20 we found a perfect solution in almost all cases. For OUR-40, the median error is high for \(n=500\) and \(n=675\). In contrast, for \(n\geq 850\), we again obtain a perfect solution in most cases. The performance of OUR-60 and OUR-80 is generally poorer, yet it shows an improving trend as the problem size increases. Across all variants, we observe that the solution improves as \(n\) increases. Conversely, the behaviour of SA differs. For small instances, SA outperforms OUR-40 to OUR-80. But the performance of SA does not improve as \(n\) increases. Consequently, all variants of the proposed algorithm perform better than SA for large problem sizes. This behaviour can be attributed to the characteristic of our decomposition method. The perfect solutions of the auxiliary NPP in the merging step are more likely to be conserved with the larger problems than the smaller ones, as depicted in \Cref{sec:theory}.

As we observe that OUR-40 already performs satisfyingly, in the following, we report the results for OUR-40 (simply referred to as OUR).

\begin{figure}[hbt!]
    \centering
    \includegraphics[width=\textwidth]{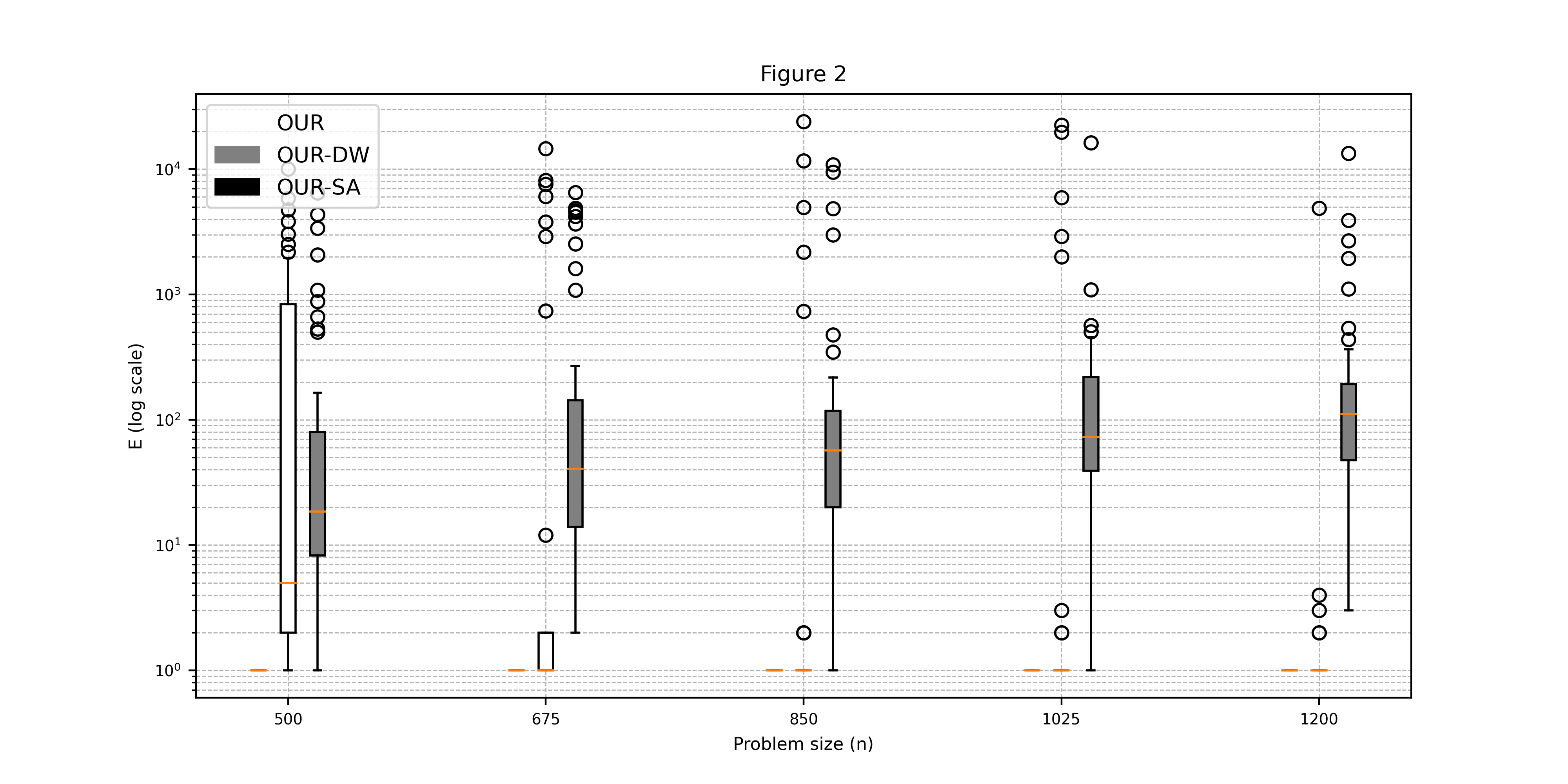}
    \caption{Experiment 2: comparing solution error \(E\) with different choice of solvers for 5 different problem sizes}
    \label{fig:2}
\end{figure}

The other adjustable parameter in \Cref{alg:partition} is the choice between QA and SA. We introduced \Cref{alg:partition} with the hybrid variant of using both QA and SA. Now we investigate this choice numerically in \Cref{fig:2} by comparing OUR (hybrid) with OUR-DW (all QA) and OUR-SA (all SA).

It is apparent that OUR-SA outperforms the other two variants. This superior performance is to be attributed to SA's ability to solve the sub-NPPs and the auxiliary NPP in the merging step with greater accuracy. First of all, an advantage arises since better solutions of sub-NPPs lead to smaller (and thus closer) values in the auxiliary NPP. The quality of the solution of the NPP in the merging step directly influences the overall solution. However, particularly for large instances, the approach can tolerate sub-optimal solutions in the sub-NPPs.  

\begin{figure}[hbt!]
    \centering
    \includegraphics[width=\textwidth]{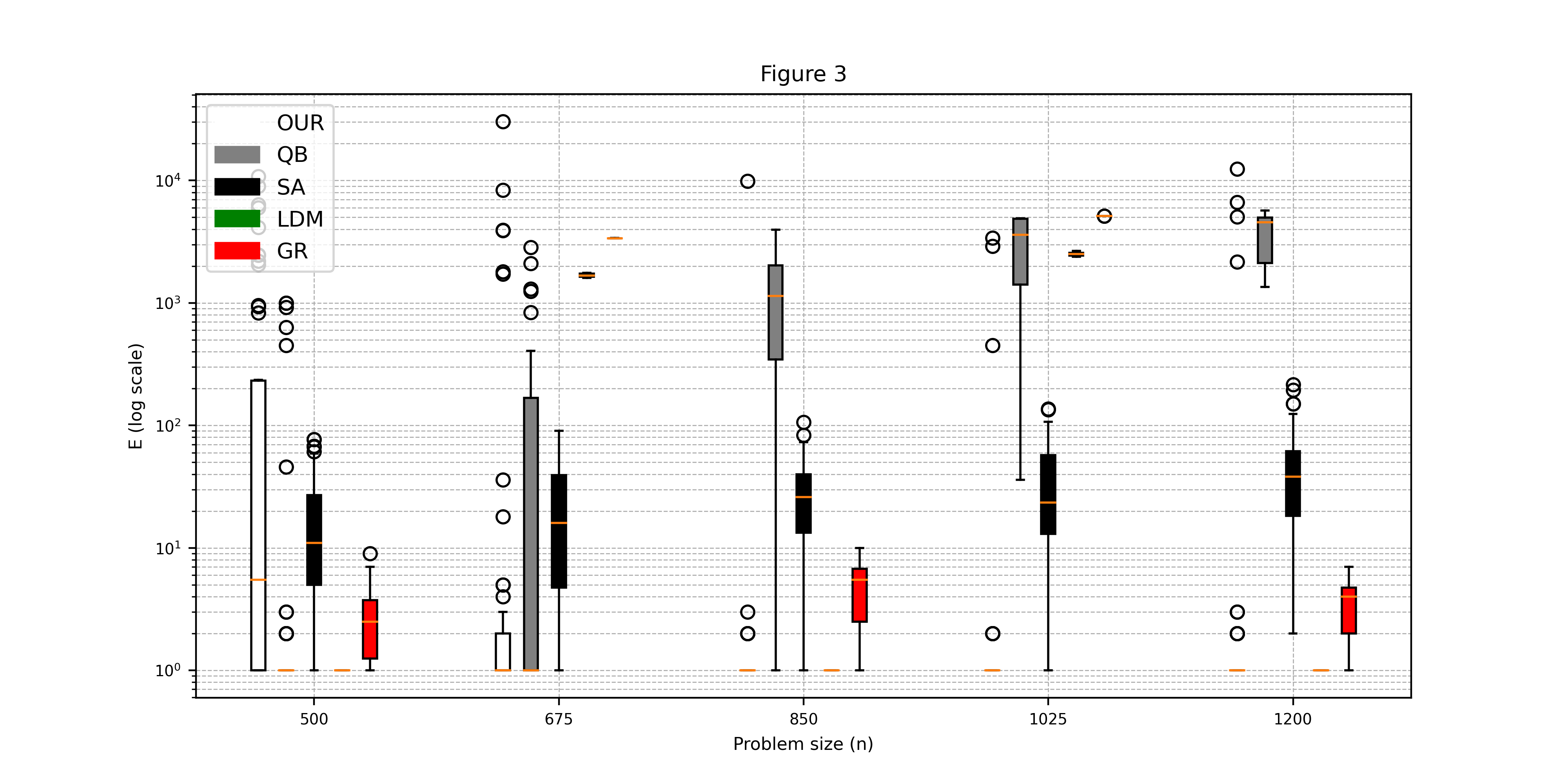}
    \caption{Experiment 3: comparing solution error \(E\) between different algorithms for 5 different problem sizes}
    \label{fig:3}
\end{figure}

In a final experiment, we compare the proposed algorithm OUR with the general decomposition solver QB. Additionally, we also evaluate OUR against the classical approaches SA, LDM and GR. The results are reported in \Cref{fig:3}. 

While QB still finds acceptable solutions for the problem sizes of \(n=500\) and \(n=675\), its solution quality declines as \(n\) increases, indicating that our method outperforms QB for large instances. Additionally, both SA and GR consistently perform worse than OUR, where GR performs better than SA for instances with even problem size. Although LDM always finds perfect solutions for every instance with an even problem size, it performs poorly for others.

For problems with an odd size, LDM and GR perform significantly worse than the hybrid OUR in \Cref{fig:3}. Furthermore, certain problems, such as \(W=\{3000, 3000, 2000,\) \( 2000, 2000\}\), are known to be ill-behaved for LDM as studied in \myCite{fisch}.

Our results suggest that \Cref{alg:partition} is promising for NPPs where the magnitude of elements is linearly related to the problem size.

Nevertheless, we note that there is plenty of room for improvements in these results. For instance, employing better embeddings specifically tailored for the NPP, such as the one introduced in \myCite{lucNPP}, might enhance the results from this study. In contrast, we use the standard \textit{EmbeddingComposite} provided by D-Wave with the default settings. Additionally, fine-tuning other parameters could also be beneficial, as demonstrated by \myCite{asproni} in their work with \textit{qbsolv}.

Besides comparing \(E\), we also investigate the potential time advantage in this study for completeness. \Cref{tab:1} shows the average time of solving a NPP of a certain size \(n\) by each algorithm.

In \Cref{tab:1}, we distinguish between the normal time and the parallel (para) time for OUR. Given that the sub-problems from \Cref{alg:partition} can be computed in parallel, we also present the time if executed in parallel. The QPU times, referring to the durations of execution on the quantum processing unit, are included in the total times.

\begin{table}[hbt!]
    \centering
    \begin{tabular}{|c|c|
    S[table-format=4.4]|
    S[table-format=4.4]|
    S[table-format=4.4]|
    S[table-format=3.4]|
    S[table-format=3.4]|}
        \hline
        \multicolumn{2}{|c|}{Algo \textbackslash \(n\)} & \textbf{500} & \textbf{675} & \textbf{850} & \textbf{1025} & \textbf{1200}\\
        \hline
        \multirow{4}{*}{\textbf{OUR-20}} & QPU & 0.7134 & 0.9473 & 1.2016 & 1.4602 & 1.7219 \\
        & total & 164.4301 & 203.8578 & 251.3105 & 339.8384 & 348.8300 \\
        & QPU (para) & 0.0347 & 0.0349 & 0.0352 & 0.0351 & 0.0352 \\
        & total (para) & 14.4464 & 13.5787 & 13.9131 & 14.4074 & 14.4290 \\
        \hline
        \multirow{4}{*}{\textbf{OUR-40}} & QPU & 0.3726 & 0.4964 & 0.6431 & 0.7729 & 0.9201 \\
        & total & 170.0892 & 219.1972 & 254.5012 & 307.5980 & 388.7510 \\
        & QPU (para) & 0.0362 & 0.0364 & 0.0364 & 0.0368 & 0.0366 \\
        & total (para) & 24.7523 & 24.7274 & 27.9884 & 37.1839 & 32.0619 \\ 
        \hline
        \multirow{4}{*}{\textbf{OUR-60}} & QPU & 0.2665 & 0.3631 & 0.4641 & 0.5665 & 0.6625 \\
        & total & 308.8090 & 443.9922 & 716.3885 & 648.3096 & 746.7421 \\
        & QPU (para) & 0.0376 & 0.0380 & 0.0383 & 0.0384 & 0.0384 \\
        & total (para) & 65.2751 & 68.5246 & 97.0904 & 62.3095 & 65.3117 \\
        \hline
        \multirow{4}{*}{\textbf{OUR-80}} & QPU & 0.2201 & 0.2936 & 0.3679 & 0.4430 & 0.5418 \\
        & total & 3060.0275 & 2943.8457 & 3311.7283 & 1573.1966 & 1306.6188 \\
        & QPU (para) & 0.0397 & 0.0402 & 0.0407 & 0.0409 & 0.0407 \\
        & total (para) & 808.9146 & 609.9785 & 551.9513 & 233.7469 & 156.1363 \\
        \hline  
        \multicolumn{2}{|c|}{\textbf{QB}} & 22.3979 & 21.1308 & 21.9535 & 25.7927 & 34.9631 \\
        \hline
        \multicolumn{2}{|c|}{\textbf{SA}} & 36.1492 & 65.2469 & 112.2975 & 175.9544 & 300.2657 \\
        \hline
        \multicolumn{2}{|c|}{\textbf{LDM}} & 0.0014 & 0.0022 & 0.0025 & 0.0035 & 0.0040 \\
        \hline
        \multicolumn{2}{|c|}{\textbf{GR}} & 0.0001 & 0.0001 & 0.0002 & 0.0002 & 0.0002 \\
        \hline
    \end{tabular}
    \caption{Rounded average time (seconds)  of 10 different problems of size \(n\) (OUR and SA are run 5 times for each problem) with \(w\in W\) drawn uniformly from \([5n,10n]\) at random}
    \label{tab:1}
\end{table}

We note here that \Cref{tab:1} is intended to illustrate the orders of magnitude of time costs, not the precise values. These times depend on many factors, which are difficult to standardize. For example, we do not implement the re-use of embeddings in OUR, although this could potentially save a huge amount of time. Additionally, waiting for quantum hardware is also included in the total times.

Nonetheless, we observe that GR always has the shortest execution time, with LDM coming in second. Even if we ignore the time necessary for embedding etc., and only compare with the QPU times, the classical heuristics still beat OUR by orders of magnitude in the present cases. However, we find that OUR-20 and OUR-40 have some time advantages over SA and, when executed in parallel, over QB as well.

\section{Conclusions}\label{sec:conclusions}
In this study, we proposed a quantum-classical algorithm to solve NPPs of a large size. We investigated its performance compared to alternative classic and hybrid algorithms. Based on numerical experiments, we found that our proposed algorithm shows better solution quality than the hybrid and general solver \textit{qbsolv} from D-Wave. This demonstrates the benefits of a special-purpose decomposition method over a general-purpose method.

Compared to classical heuristics, our algorithm is competitive in terms of solution quality. Especially in instances with an odd problem size, our proposed method finds significantly better solutions. A time advantage of our algorithm is not observable from the results in this study. Although our implementation can be improved to achieve shorter times, even the quantum annealing time in the case of parallel execution is longer than the execution time of the classical heuristics.

We note that this study focuses on the 2-way NPP. Nonetheless, by introducing additional quadratic constraints to the auxiliary NPP in the last step of our algorithm (quadratic constraints can be included in a QUBO as described in \myCite[p. 117]{alidaee}), our method could be generalized to multi-way NPPs. Another generalization can be made in terms of problem size. For example, multiple layers of sub-problems can be introduced. Furthermore, this study only examined an uniformly random decomposing vector \(j\). It might be beneficial to investigate other \(j\), such as one which considers the magnitude of the elements. Moreover, the performance of our proposed algorithm for other NP-hard problems can also be investigated by reducing these problems to NPPs.

As a further alternative, different solvers can be employed within the method instead of QA and SA. For example, as an alternative to QA, QUBOs can also be solved on gate-based quantum devices using QAOA as described, for example, in \myCite{Hadfield2019}. 

Lastly, we note that future improvements of quantum hardware will most likely enhance the performance of our method, both with respect to execution time and quality of solutions, without any modification of our proposed method.

\section{Acknowledgments}\label{sec:acknow}
We would like to thank Dominik Leib for valuable discussions and help. This work was partially funded by the German Federal Ministry of Education and Research (Bundesministerium für Bildung und Forschung) within the project \textit{Rymax One}, and by \textit{The Quantum Initiative Rhineland-Palatinate QUIP}. Numerical experiments were carried out on the \textit{ETH Euler Cluster} (Swiss Federal Institute of Technology Zurich) by Zongji Li as a guest user.

%% file: main.bib
@article{kochenberger2014,
        title={The unconstrained binary quadratic programming problem: a survey},
        author={Kochenberger, Gary and Hao, Jin-Kao and Glover, Fred and Lewis, Mark and L{\"u}, Zhipeng and Wang, Haibo and Wang, Yang},
        journal={Journal of Combinatorial Optimization},
        year={2014},
        doi={10.1007/s10878-014-9734-0},
}

@misc{dwave2020,
	author = {{D-Wave Systems Inc.}},
	year   = {2020},
	title  = {{D-Wave} Hybrid Solver Service: An Overview},
	url = {https://www.dwavesys.com/media/4bnpi53x/14-1039a-b_d-wave_hybrid_solver_service_an_overview.pdf}, 
	note   = {Accessed 2022-11-08}
}

@article{alexeev2021,
        title = {Quantum Computer Systems for Scientific Discovery},
        author = {Alexeev, Yuri and Bacon, Dave and Brown, Kenneth R. and Calderbank, Robert and Carr, Lincoln D. and Chong, Frederic T. and DeMarco, Brian and Englund, Dirk and Farhi, Edward and Fefferman, Bill and Gorshkov, Alexey V. and Houck, Andrew and Kim, Jungsang and Kimmel, Shelby and Lange, Michael and Lloyd, Seth and Lukin, Mikhail D. and Maslov, Dmitri and Maunz, Peter and Monroe, Christopher and Preskill, John and Roetteler, Martin and Savage, Martin J. and Thompson, Jeff},
        journal = {PRX Quantum},
        year = {2021},
        doi = {10.1103/PRXQuantum.2.017001},
}

@article{krasecki2023,
        title={The Role of Experimental Noise in a Hybrid Classical-Molecular Computer to Solve Combinatorial Optimization Problems},
        author={Krasecki, Veronica K. and Sharma, Abhishek and Cavell, Andrew C. and Forman, Christopher and Guo, Si Yue and Jensen, Evan Thomas and Smith, Mackinsey A. and Czerwinski, Rachel and Friederich, Pascal and Hickman, Riley J. and Gianneschi, Nathan and Aspuru-Guzik, Alán and Cronin, Leroy and Goldsmith, Randall H.},
        journal={ACS Central Science},
        year={2023},
        doi={10.1021/acscentsci.3c00515},
}

@article{asproni,
        title = {Accuracy and minor embedding in subqubo decomposition with fully connected large problems: a case study about the number partitioning problem},
        author = {Asproni, Luca and Caputo, Davide and Silva, Blanca and Fazzi, Giovanni and Magagnini, Marco},
        journal = {Quantum Machine Intelligence},
        year = {2020},
        doi = {10.1007/s42484-020-00014-w}
}

@article{lucNPP,
        author = {Lucas, Andrew},
        title = {Hard combinatorial problems and minor embeddings on lattice graphs},
        journal = {Quantum Information Processing},
        year = {2019},
        doi = {10.1007/s11128-019-2323-5}
}

@InProceedings{zbinden,
        author = {Zbinden, Stefanie
        and Bärtschi, Andreas
        and Djidjev, Hristo
        and Eidenbenz, Stephan},
        title = {Embedding Algorithms for Quantum Annealers with Chimera and Pegasus Connection Topologies},
        booktitle = {International Conference on High Performance Computing},
        year = {2020},
        doi = {10.1007/978-3-030-50743-5_10}
}

@article{bass,
        author = {Bass, Gideon and Henderson, Maxwell and Heath, Joshua and Dulny, Joseph},
        title = {Optimizing the optimizer: decomposition techniques for quantum annealing},
        journal = {Quantum Machine Intelligence},
        year = {2021},
        doi = {10.1007/s42484-021-00039-9}
}

@article{korf98,
        title = {A complete anytime algorithm for number partitioning},
        journal = {Artificial Intelligence},
        year = {1998},
        doi = {10.1016/S0004-3702(98)00086-1},
        author = {Richard E. Korf}
}

@article{schroeppel,
        author = {Schroeppel, Richard and Shamir, Adi},
        title = {A \(T = O(2^{n/2} )\), \(S = O(2^{n/4} )\) Algorithm for Certain NP-Complete Problems},
        journal = {SIAM Journal on Computing},
        year = {1981},
        doi = {10.1137/0210033}
}

@book{Pinedo.2012,
         author = {Pinedo, Michael L.},
         year = {2012},
         title = {Scheduling: Theory, Algorithms, and Systems},
         doi = {10.1007/978-1-4614-2361-4}
}

@article{horowitz,
        author = {Horowitz, Ellis and Sahni, Sartaj},
        title = {Computing Partitions with Applications to the Knapsack Problem},
        year = {1974},
        doi = {10.1145/321812.321823},
        journal = {Journal of the ACM}
}

@techreport{kk,
        author = {Karmarkar, Narendra and Karp, Richard M.},
        title = {The Differencing Method of Set Partitioning},
        institution = {EECS Department, University of California, Berkeley},
        year = {1983},
        uRL = {http://www2.eecs.berkeley.edu/Pubs/TechRpts/1983/6353.html},
        number = {UCB/CSD-83-113}
}

@misc{sa,
        author = {{D-Wave Systems Inc.}},
        title = {Simulated Annealing Sampler - dwave-neal 0.5.9 documentation},
        url = {https://docs.ocean.dwavesys.com/projects/neal/en/latest/reference/sampler.html},
        note = {Accessed on 2023-09-16},
        year = {2023}
}

@article{analysis,
        author = {Gent, Ian P. and Walsh, Toby},
        title = {Analysis of Heuristics for Number Partitioning},
        journal = {Computational Intelligence},
        doi = {10.1111/0824-7935.00069},
        year = {1998}
}

@article{alidaee,
        author = {Alidaee, Bahram and Glover, Fred and Kochenberger, Gary A. and Rego, Cesar},
        title = {A new modeling and solution approach for the number partitioning problem},
        journal = {Advances in Decision Sciences},
        doi = {10.1155/JAMDS.2005.113},
        year = {2005}
}

@article{coffman,
        author = {Coffman, E. G. and Johnson, D. S. and Lueker, G. S. and Shor, P. W.},
        journal = {Statistical Science},
        title = {Probabilistic Analysis of Packing and Related Partitioning Problems},
        doi = {10.1214/ss/1177011082},
        year = {1993}
}

@article{boettcher,
        author = {Boettcher, S. and Mertens, S.},
        title = {Analysis of the Karmarkar-Karp differencing algorithm},
        journal = {The European Physical Journal B},
        doi = {10.1140/epjb/e2008-00320-9},
        year = {2008}
}

@article{michiels,
        author = {Michiels, Wil and Korst, Jan and Aarts, Emile and van Leeuwen, Jan},
        title = {Performance ratios of the Karmarkar-Karp differencing method},
        journal = {Journal of Combinatorial Optimization},
        doi = {10.1007/s10878-006-9010-z},
        year = {2007}
}

@incollection{mertens,
        author = {Mertens, Stephan},
        title = "{The Easiest Hard Problem: Number Partitioning}",
        booktitle = "{Computational Complexity and Statistical Physics}",
        year = {2005},
        doi = {10.1093/oso/9780195177374.003.0012}
}

@article{qubom,
        author = {Glover, Fred and Kochenberger, Gary and Hennig, Rick and Du, Yu},
        title = {Quantum bridge analytics I: a tutorial on formulating and using QUBO models},
        journal = {Annals of Operations Research},
        year = {2022},
        doi = {10.1007/s10479-022-04634-2}
}

@article{fisch,
        author = {Fischetti, Matteo and Martello, Silvano},
        title = {Worst-case analysis of the differencing method for the partition problem},
        journal = {Mathematical Programming},
        year = {1987},
        doi = {10.1007/BF02591687}
}

@article{okada,
        author = {Okada, Shuntaro and Ohzeki, Masayuki and Terabe, Masayoshi and Taguchi, Shinichiro},
        title = {Improving solutions by embedding larger subproblems in a D-Wave quantum annealer},
        journal = {Scientific Reports},
        year = {2019},
        doi = {10.1038/s41598-018-38388-4}
}

@inproceedings{osaba,
        author={Osaba, Eneko and Villar-Rodriguez, Esther and Oregi, Izaskun and Moreno-Fernandez-de-Leceta, Aitor},
        booktitle={IEEE Congress on Evolutionary Computation}, 
        title={Hybrid Quantum Computing - Tabu Search Algorithm for Partitioning Problems: Preliminary Study on the Traveling Salesman Problem}, 
        year={2021},
        doi={10.1109/CEC45853.2021.9504923}
}

@article{raymond23,
        author = {Raymond, Jack and Stevanovic, Radomir and Bernoudy, William and Boothby, Kelly and McGeoch, Catherine C. and Berkley, Andrew J. and Farr\'{e}, Pau and Pasvolsky, Joel and King, Andrew D.},
        title = {Hybrid Quantum Annealing for Larger-than-QPU Lattice-Structured Problems},
        year = {2023},
        doi = {10.1145/3579368},
        journal = {ACM Transactions on Quantum Computing},
}

@article{shay,
        author={Shaydulin, Ruslan and Ushijima-Mwesigwa, Hayato and Negre, Christian F. A. and Safro, Ilya and Mniszewski, Susan M. and Alexeev, Yuri},
        journal={Computer}, 
        title={A Hybrid Approach for Solving Optimization Problems on Small Quantum Computers}, 
        year={2019},
        doi={10.1109/MC.2019.2908942}
}

@article{atobe,
        author={Atobe, Yuta and Tawada, Masashi and Togawa, Nozomu},
        journal={IEEE Transactions on Computers}, 
        title={Hybrid Annealing Method Based on subQUBO Model Extraction With Multiple Solution Instances}, 
        year={2022},
        doi={10.1109/TC.2021.3138629}
}

@techreport{dwave1,
        author = {Booth, Michael and Reinhardt, Steven P. and Roy, Aidan},
        title = {Partitioning Optimization Problems for Hybrid Classical/Quantum Execution},
        institution = {D-Wave Systems Inc.},
        year = {2017},
        url = {https://www.dwavesys.com/media/jhlpvult/partitioning_qubos_for_quantum_acceleration-2.pdf},
        number = {14-1006A-A}
}

@inproceedings{npcomplete,
	author = {Karp, Richard M.},
	title = {Reducibility among Combinatorial Problems},
	booktitle = {Complexity of Computer Computations},
	doi={10.1007/978-1-4684-2001-2_9},
	year = {1972}
}

@misc{abbas2023quantum,
      title={Quantum Optimization: Potential, Challenges, and the Path Forward}, 
      author={Amira Abbas and Andris Ambainis and Brandon Augustino and Andreas Bärtschi and Harry Buhrman and Carleton Coffrin and Giorgio Cortiana and Vedran Dunjko and Daniel J. Egger and Bruce G. Elmegreen and Nicola Franco and Filippo Fratini and Bryce Fuller and Julien Gacon and Constantin Gonciulea and Sander Gribling and Swati Gupta and Stuart Hadfield and Raoul Heese and Gerhard Kircher and Thomas Kleinert and Thorsten Koch and Georgios Korpas and Steve Lenk and Jakub Marecek and Vanio Markov and Guglielmo Mazzola and Stefano Mensa and Naeimeh Mohseni and Giacomo Nannicini and Corey O'Meara and Elena Peña Tapia and Sebastian Pokutta and Manuel Proissl and Patrick Rebentrost and Emre Sahin and Benjamin C. B. Symons and Sabine Tornow and Victor Valls and Stefan Woerner and Mira L. Wolf-Bauwens and Jon Yard and Sheir Yarkoni and Dirk Zechiel and Sergiy Zhuk and Christa Zoufal},
      year={2023},
      eprint={2312.02279},
      archivePrefix={arXiv},
      primaryClass={quant-ph}
}

@article{Hadfield2019,
   title={From the Quantum Approximate Optimization Algorithm to a Quantum Alternating Operator Ansatz},
   volume={12},
   DOI={10.3390/a12020034},
   journal={Algorithms},
   author={Hadfield, Stuart and Wang, Zhihui and O’Gorman, Bryan and Rieffel, Eleanor and Venturelli, Davide and Biswas, Rupak},
   year={2019}
}
